\documentclass[10pt]{amsart}
\usepackage{amsfonts,amssymb,amsmath,amscd,amsthm}
\usepackage[all]{xy}

\providecommand{\tmap}[2]{\ensuremath{#1_{#2}}}
\providecommand{\cprod}[2]{\ensuremath{ #1\rtimes_{#2}\mathbb{Z}}}
\providecommand{\ds}[1]{{\displaystyle{#1}}}
\providecommand{\ts}[1]{{\textstyle{#1}}}

\newcommand{\C}{{\mathbb C}}

\newcommand{\Z}{{\mathbb Z}}
\newcommand{\R}{\mathbb R}
\newcommand{\Q}{{\mathbb Q}}

\newcommand{\ttau}{\overline{\tau}}

\newcommand{\id}{\mathrm{Id}}
\newcommand{\im}{\mathrm{Im}}
\newcommand{\ev}{\mathrm{ev}}
\newcommand{\Aut}{\mathrm{Aut}}
\renewcommand{\d}{\mathrm{d}}
\newcommand{\entr}{\mathrm{ht}}
\newcommand{\topentr}{\mathrm{h}_{\mathrm{top}}}

\newcommand{\calH}{\mathcal{H}}
\newcommand{{\frakM}}{\mathfrak{M}}
\newcommand{{\frakS}}{\mathfrak{S}}
\newcommand{\calBH}{\mathcal{B}(\mathcal{H})}

\renewcommand{\phi}{\varphi}

\newtheorem{theorem}{Theorem}[section]
\newtheorem{lemma}[theorem]{Lemma}
\newtheorem{corollary}[theorem]{Corollary}
\newtheorem{proposition}[theorem]{Proposition}
\newtheorem{definition}[theorem]{Definition}
\newtheorem{example}[theorem]{Example}

\title[The C*- algebras Associated to Time-$t$ Automorphisms]{The C*-algebras Associated to Time-$t$ Automorphisms of  Mapping Tori}
\author{Benjam\'{\i}n Itz\'a-Ortiz} 
\date{February 8th,~2005}

\begin{document}

\begin{abstract}
We find the range of a trace on the  $K_0$ group of a  crossed product by a   time-$t$ automorphism of a mapping torus.  We also find a formula to compute the Voiculescu-Brown entropy for such an automorphism. By specializing to the commutative setting, we prove that the crossed products by minimal time-t homeomorphisms of suspensions built over strongly orbit equivalent Cantor minimal systems have isomorphic Elliott invariants.   As an application of our results we give examples of  dynamical systems on (compact metric) connected 1-dimensional spaces
which are not flip conjugate (because of different entropy) yet their associated crossed products have isomorphic Elliott invariants. 
\end{abstract}

\maketitle


\section*{Introduction} 
 Let $A$ be a   C*-algebra and let $\alpha$ be an automorphism of $A$. The {\em mapping torus} of $\alpha$ (see, e.g. \cite[V.10.3]{Black}) is the  C*-algebra
 \[
    M\sb\alpha=\{ f\colon [0,1]\rightarrow A \, |\, f(1)=\alpha(f(0))\}.
 \] 
There is a natural action ${\phi}$ of $\R$ on $M_\alpha$ given by translation. That is, 
for $t\in \R$, $f\in M\sb\alpha$ and $s\in [0,1]$,
\[
     (\phi_t(f))(s) = \alpha^n(f(s+t-n)) 
\]
where $n$ is the unique integer such that $n\leq s+t <n+1$. We call $(M\sb\alpha,{\phi})$ the {\em mapping torus flow} associated to $(A,\alpha)$. For a fixed $t\in \R$ we call $\phi_t$ the time-$t$ automorphism of $M\sb\alpha$.  
 
The mapping torus flow is an extension to the non-commutative setting of the suspension flow, a standard construction in topological dynamics which we describe next, cf. \cite[II.5.5]{V}.
We say that $(X,S)$ is a {\em dynamical system} if $X$ is a compact metric space and $S$ is a homeomorphism of $X$.   We call the {\em suspension} of $S$  the quotient space $Y$ obtained from the space $X\times [0,1]$ by identifying points
$(x,1)$ and $(S(x),0)$ with each other for all $x\in X$. It turns out that $Y$ is also a compact metric space and has (covering) dimension $\dim(Y)=\dim(X)+1$. Let $\pi\colon X\times [0,1]\rightarrow Y$ be the canonical quotient map and define $\pi(x,s)=[x,s]$ for $(x,s)\in X\times [0,1]$.  There is a natural action $T$ of $\R$ on $Y$ given by translation along the second coordinate. That is, for $t\in \R$ and  $[x,s]\in Y$, 
\[
     T^t\left([x,s]\right)=\left[S^n(x), s+t-n\right]
\]
where $n$ is the unique integer such that $n\leq s+t <n+1$. We call 
$(Y,T)$ the {\em suspension flow} associated to  $(X,S)$.  We sometimes will refer to $(X,S)$ (or just $S$) as the base of the suspension $(Y,T)$ and to $Y$ as the  suspension built over $(X,S)$ (or just $S$).  For a fixed $t\in\R$ we call $T^t$ the time-$t$ homeomorphism of $Y$.

If $A$ is a C*-algebra and $\alpha$ is an automorphism of $A$ then the  {\em crossed product} $\cprod{A}{\alpha}$ induced by $\alpha$ is the (universal)  C*-algebra generated by $A$ and  a unitary $u$ satisfying $u^nau^{-n}=\alpha^n(a)$ for all $a\in A$ and for all $n\in\Z$. 

Recall that a dynamical system $(X,S)$ induces a  C*-dynamical system $(C(X)$, $\alpha)$ where $\alpha$ is the automorphism of $C(X)$ defined by $\alpha(f)=f\circ S$.   We define 
$C\sp\ast(\Z,X,S)$ to be the crossed product $\cprod{C(X)}{\alpha}$.
If $T$ is an action of $\R$ on a compact metric space  $Y$ then $(Y,T)$ induces a C*-dynamical $\R$-system $(C(Y),\phi)$  where ${\phi}$ is the homomorphism from $\R$ to $\Aut(C(Y))$ given by $\phi_t(g)=g\circ T^t$ for $t\in\R$ and $g\in C(Y)$. Suppose that  $(Y,T)$ is the suspension flow of $(X,S)$. Then, if
$(M\sb\alpha,{\phi})$ is the mapping torus flow associated to  $(C(X),\alpha)$,  it turns out that, for each $t\in T$, the C*-algebras  $C\sp\ast (\Z,Y,T^t)$ and $\cprod{M\sb\alpha}{\phi_t}$ are isomorphic. Therefore, by studying the crossed products by time-$t$ automorphisms of mapping tori, we are studying the extension to the non-commutative setting of crossed products by time-$t$ homeomorphisms of suspension spaces.

It is well known that two minimal homeomorphisms of the unit circle $S^1$ are flip conjugate if and only if their associated crossed products are isomorphic. 
Assuming the Elliott conjecture, which says that a complete isomorphism invariant for the associated crossed products is of $K$-theoretic nature, we find examples showing that this good behavior on $S^1$ does not extend to other compact metric connected 1-dimensional spaces; namely, we find examples of minimal homeomorphisms on compact metric connected 1-dimensional spaces which are not flip conjugate yet their associated crossed products have isomorphic Elliott invariants. Since (topological) orbit equivalence implies flip conjugacy for minimal homeomorphisms of compact metric connected spaces, our examples are not even orbit equivalent. Therefore, at the same time, we are showing that  the Krieger type theorem  of Giordano, Putnam and Skau \cite[Theorem~2.1]{GPS} for the family of minimal Cantor systems does not extend to a larger class of minimal dynamical systems. Although this had been established, we remark that the  known examples  involve spaces which are either disconnected or have dimension at least 2 (see \cite{NCP1} for a survey).


We divide this paper in 3 sections. In Section~1 we find the range of a trace on the 
$K_0$ group of a crossed product by a time-$t$ automorphism of a mapping torus.   In Section~2 we will give a formula to compute the Voiculescu-Brown entropy for such  an automorphism.  In Section~3 we specialize our results to the commutative setting. In particular, we show that crossed products by minimal time-$t$ homeomorphisms of suspensions built over strongly orbit equivalent Cantor minimal systems have isomorphic  Elliott invariants. When the base of the suspension is minimal we prove that, except for a countable set containing $\Q$, all other time-$t$ homeomorphisms are minimal.  As applications of our results we provide  examples showing  the  existence of dynamical systems on (compact metric) connected 1-dimensional spaces which are not flip conjugate (because of different entropy) yet their associated  crossed products have isomorphic Elliott invariants.

The author wishes to express his gratitude to his Ph.D.\ advisor, Professor N.\ Christopher Phillips, for suggesting this line of research. Above all, we are grateful to him for his continuous encouragement and for being a constant source of inspiration. We would also like to thank Professors Thierry Giordano, David Handelman and Vladimir Pestov for their support, encouragement and for many stimulating conversations.
\section{Ordered $K_0$-group}
The order structure on the $K_0$ group of a C*-algebra has played an important role in many of the applications of $K$-theory to C*-algebras, particularly to the structure of AF algebras (see eg.\ \cite[III.7]{Black}). It is usually a nontrivial problem to determine the order structure on the $K_0$ group of a crossed product $\cprod{A}{\alpha}$. However, by calculating the range of the states on $K_0(\cprod{A}{\alpha})$ which come from traces on the crossed product, one at least may recover partial information about the order  structure of $K_0(\cprod{A}{\alpha})$.  In this section we will give a formula to compute the range of a trace on the $K_0$ group of a crossed product by a time-$t$ automorphism of a mapping torus.

When dealing with the $K$-theory of crossed products, one has two fundamental tools (see eg.\ \cite[V.10.2]{Black}), namely:
\begin{theorem}[Pimsner-Voiculescu exact sequence]
Let $A$ be a C*-algebra and let $\alpha$ be an automorphism of $A$. Then there is a natural cyclic six-term exact sequence
\[
\xymatrix {%
      K_0(A) \ar[r]^{\id-\alpha\sb\ast} & K_0(A) \ar[r]^(.43){i\sb\ast} & K_0(\cprod{A}{\alpha})
                                                                                          \ar[d]^{\partial}\\
      K_1(\cprod{A}{\alpha}) \ar[u]    & K_1(A)\ar[l]_(.38){i\sb\ast}   & K_1(A)
                                                                                          \ar[l]_(.43){\id-\alpha\sb\ast}\\
               }
\]
where the maps $i\sb\ast\colon K_j(A)\rightarrow K_j(\cprod{A}{\alpha})$ for $j=0,1$,  are the maps on $K$-theory induced by the inclusion $A\rightarrow \cprod{A}{\alpha}$.
\end{theorem}
\begin{theorem}[Connes' Thom Isomorphism]
If $\phi\colon \R\rightarrow \mathrm{Aut}(A)$ is an action of $\R$ on the C*-algebra $A$ then $K_i(A\rtimes\sb\phi\R)\cong K_{1-i}(A)$ for $i=0,1$.
\end{theorem}

As a corollary of his  Thom Isomorphism, Connes \cite[Corollary~V.6]{C} proved the following important result.

\begin{theorem}[Connes' Isomorphism]\label{ConnesIsom}
Let $\alpha$ be an automorphism of a C*-algebra $A$ and let $M_\alpha$ be the mapping torus of $\alpha$. Then $K_i(M_\alpha)\cong K_{1-i}(\cprod{A}{\alpha})$ for $i=0,1$.
\end{theorem}

Explicit formulas for Connes'  Isomorphisms (\ref{ConnesIsom}) have been found by Pashke \cite{P}. (See also \cite[Proposition V.10.4.3]{Black} and \cite[V.4]{E}.) One of Connes' Isomorphisms will be used often in this section so it will be convenient to label it.  We do this in the following.
\begin{definition}\label{connesk}
 Let $\alpha$ be an automorphism of a C*-algebra $A$ and let $M\sb\alpha$ be the mapping torus of $\alpha$. We will denote by
\[
           k\sb\alpha\colon K_1(M\sb\alpha)\rightarrow K_0(\cprod{A}{\alpha})
\]
the corresponding Connes' Isomorphism of Theorem~\ref{ConnesIsom}.         
\end{definition}
  
We now define the following homomorphism for later reference. It relates the $K_1$-groups of a mapping torus and a ``double" mapping torus.

\begin{definition}\label{theta}
Let $A$ be a C*-algebra, let $\alpha\in\Aut(A)$ and let $(M\sb\alpha,\phi)$ be the mapping torus flow associated to $(A,\alpha)$.   Let  $t\in\R$ and let  $M_{\phi_t}$ be the mapping torus of $\phi_t$. Define a homomorphism
    \[
       \theta_t\colon K_1(M\sb\alpha)\rightarrow K_1(M_{\phi_t})
    \] 
by $\theta_t([u])=[\phi_{st}(u)]$, where $u$  is a unitary representing an element in $K_1(M\sb\alpha)$ and $\phi_{st}(u)$ stands for the map $s\mapsto \phi_{st}(u)$.
  \end{definition}
  
  It is straightforward to check that the map $\theta_t$ of Definition~\ref{theta} in fact gives a homomorphism (of abelian groups).
  
  Let $A$ be a C*-algebra, let $\alpha\in\Aut(A)$ and let  $(M\sb\alpha,{\phi})$ be  the mapping torus flow associated to $(A,\alpha)$. 
For $t\in\R$ consider the time-$t$ automorphism $\phi_t$ of $M\sb\alpha$. 
Since $\phi_t$ is homotopic to the identity, the Pimsner-Voiculescu exact sequence gives us the short exact sequence

\begin{equation}\label{ses}
\xymatrix{   0 \ar[r]& K_0(M\sb\alpha) \ar[r]^(.43){i\sb\ast}
              & K_0( \cprod{M\sb\alpha}{\tmap{{\phi}}{t}}) \ar[r]^(.57){\partial}& K_1(M\sb\alpha) \ar[r]& 0
                }
\end{equation}
where $i\sb\ast:K_0(M\sb\alpha)\rightarrow K_0(\cprod{M\sb\alpha}{\phi_t})$ is the map induced by the inclusion $M\sb\alpha\rightarrow \cprod{M\sb\alpha}{\phi_t}$.

%
%
\begin{lemma}\label{lemmasplit}
The short exact sequence (\ref{ses}) splits.  A splitting map is given by the  composition
\[
      k_{\phi_t}\circ\theta_t\colon K_1(M\sb\alpha)\rightarrow K_0(\cprod{M\sb\alpha}{\phi_t})
\]
where $\theta_t\colon K_1(M\sb\alpha)\rightarrow K_1(M_{\phi_t})$ is the map of Definition~\ref{theta} and $k_{\phi_t}\colon K_1(M_{\phi_t})\rightarrow K_0(\cprod{M\sb\alpha}{\phi_t})$ is the Connes' Isomorphism (\ref{connesk}).
\end{lemma}
\begin{proof}
Let $(\ev_0)\sb\ast \colon  K_1( M_{\phi_t} )\rightarrow K_1(M\sb\alpha)$ be the 
 map  induced by the evaluation at zero map $\ev_0\colon M_{\phi_t}\rightarrow M\sb\alpha$.  As a consequence of \cite[Lemma 1]{P}, we have that 
 $\partial\circ k_{\phi_t} =(\ev_0)\sb\ast$.  Hence, if $u$ is a unitary representing 
 an element in $K_1(M\sb\alpha)$  then
    \begin{eqnarray*}
      ( \partial \circ k_{\phi_{t}}\circ\theta_t)([u]) & = & (\ev_0)\sb\ast (\theta_t ([u])) \\
             & = & (\ev_0)\sb\ast ( [\phi_{st}(u)]) \\
             & = & [u].
    \end{eqnarray*}
\end{proof}

%
%
\begin{theorem}\label{kgroup}
Let $A$ be a C*-algebra, let $\alpha$ be an automorphism on $A$ and let $(M\sb\alpha,{\phi})$ be the mapping torus flow associated to $(A,\alpha)$. If $t\in\R$ then the  map
\[
 K_0(M\sb\alpha)\oplus K_0 (\cprod{A}{\alpha}) \xrightarrow{i\sb\ast +k_{\phi_{t}}\circ\theta_t\circ k_{\alpha}^{-1}} K_0(\cprod{ M\sb\alpha }{ \phi_t })
\]
is an isomorphism, where $i\sb\ast \colon K_0(M\sb\alpha)\rightarrow K_0(\cprod{M\sb\alpha}{\phi_t})$ is the map induced by the inclusion $M\sb\alpha\rightarrow
\cprod{ M\sb\alpha }{ \phi_t }$,
the maps $k_{\phi_{t}}\colon K_1(M_{\phi_t})\rightarrow K_0(\cprod{M\sb\alpha}{\phi_t})$ and  $k_{\alpha}\colon K_1(M\sb\alpha)\rightarrow K_0(\cprod{A}{\alpha})$  are Connes' Isomorphisms (\ref{connesk})  and $\theta_t\colon K_1(M\sb\alpha)\rightarrow K_1(M_{\phi_t})$ is as in Definition~\ref{theta}.
\end{theorem}
\begin{proof}
 From Lemma~\ref{lemmasplit}, the  map $i\sb\ast + k_{\phi_t}\circ\theta_t$ is an 
 isomorphism from $K_0(M\sb\alpha)\oplus K_1(M\sb\alpha)$ onto  $K_0(\cprod{M\sb\alpha}{\phi_t})$.
 Compose $i\sb\ast + k_{\phi_t}\circ\theta_t$ with the isomorphism $\id\oplus k_{\alpha}^{-1}\colon K_0(M\sb\alpha)\oplus K_0(\cprod{A}{\alpha})
 \rightarrow K_0(M\sb\alpha)\oplus K_1(M\sb\alpha)$, where $k_{\alpha}$ is Connes' Isomorphism (\ref{connesk}),  to obtain the
 desired result.
\end{proof}

The isomorphism of Theorem~\ref{kgroup} will be crucial in obtaining a formula to calculate the range of states on $K_0(\cprod{M\sb\alpha}{\phi_t})$  which come from traces on the crossed product $\cprod{A}{\alpha}$.
In Section~3 we will present some good cases for which our results in this section yield complete information about the order structure of $\cprod{M\sb\alpha}{\phi_t}$.

For the rest of this section we will consider only unital C*-algebras, so as to ensure that   the order structure and states on  $K_0$ groups make sense. 
Recall that a trace $\tau$ on a unital  C*-algebra $A$ is a linear function $\tau\colon A\rightarrow \C$ satisfying $\tau(ab)=\tau(ba)$ and $\tau(a\sp\ast a)\geq 0$
for all $a$ and $b$ in $A$.  We also require that $\tau (1)=1$, that is,  our traces are normalized.  For all $n\geq 1$,  a trace $\tau$  extends naturally to the *-algebra $M_n(A)$ of all $n\times n$ matrices over $A$ by the formula $\tau( a )=\sum_{i=1}^{n}\tau(a_{i,i})$ for all
$a=(a_{i,j})_{1\leq i,j\leq n}$ in $M_n(A)$. A trace $\tau$ then defines a group  homomorphims $\tau\colon K_0(A)\rightarrow \R$, via the formula $\tau ([p] - [q]) = \tau (p) - \tau (q)$. 
Given a trace preserving automorphism $\alpha$ of $A$,  $\tau$ extends canonically to a trace on $\cprod{A}{\alpha}$ by the formula $\tau(\sum a_n u^n)=\tau (a_0)$.  It turns out that any other extension of $\tau $ to $\cprod{A}{\alpha}$ defines the same state on $K_0(\cprod{A}{\alpha})$. Conversely, any trace on $\cprod{A}{\alpha}$  restricts to an $\alpha$-invariant trace on $A$. Observe that we are using the same notation $\tau$ for the original trace, its natural extensions and the corresponding state on $K_0$.

\begin{definition}\label{ttau}
Let $A$ be a unital C*-algebra, let $\alpha$ be an automorphism of $A$  and let $(M\sb\alpha,{\phi})$ be the mapping torus flow associated to $(A,\alpha)$. If  $\tau$ is an $\alpha$-invariant trace on $A$ then $\tau$ induces a ${\phi}$-invariant trace $\ttau$ on $M\sb\alpha$ via the formula
\[
   \ttau(f)=\int_{[0,1]}\tau(f(s)) \, \mathrm{d}\lambda(s).
\]
for $f\in M\sb\alpha$ and where $\lambda$ is  Lebesgue measure on $[0,1]$.
\end{definition}

It is straightforward to check that the map $\ttau$ of Definition~\ref{ttau} is in fact a $\phi$-invariant trace.
We now prove a technical lemma. Recall that $U\sb\infty(A)$ denotes the unitary group of $A$  defined to be the inductive limit of the sequence of groups
\[
   \xymatrix{ U_1(A)  \ar[r]^{i_1}& U_2(A)  \ar[r]^(.58){i_2}&  \cdots  \ar[r]^(.43){i_{n-1}} &
                   U_n(A)  \ar[r]^(.45){i_n}    & U_{{n+1}}(A)\ar[r]^(.6){i_{n+1}} &\cdots
                 }
\]
where, for all $n\geq 1$, $U_n(A)$ is the subgroup of unitary matrices of $M_n(A)$ and 
$i_n$ maps $a$ in $U_n(A)$ to $\bigl( \begin{smallmatrix} a & 0 \\ 0 & 1 \end{smallmatrix}\bigr)$ in $U_{n+1}(A)$.

\begin{lemma}\label{derivativeofphi}
 Let  the notation be as in Definition~\ref{ttau}.  Let $u$ be an element in  $U\sb\infty (M\sb\alpha)$.
  Viewing $u$ as a function $u\colon [0,1]\to U\sb\infty(A)$, assume that $u$ is smooth and that $u$ is constant on small neighborhoods of $0$ and $1$. Let $t\in\R$ and let $M_{\phi_t}$ be the mapping torus of $\phi_t$. Denote by $\phi_{st}(u)\in M_{\phi_t}$  the map defined by $s\mapsto \phi_{st}(u)$. Then $\phi_{st}$ is smooth and 
   \[
      \frac{\d}{\d s}\phi_{st} (u) = t\, \phi_{st} (u\sp\prime ),
   \]
where $t\, \phi_{st}(u\sp\prime)$ stands for the map $s\mapsto t\, \phi_{st}(u\sp\prime)$.
\end{lemma}
\begin{proof}
     For $0\leq r\leq 1$, by definition we have
 \[
      \phi_{st}(u)(r)=\alpha^n\left(u(r+st-n)\right),
 \]
where $n$ is the unique integer such that $0\leq r+st-n <1$. The lemma will then follow after a standard computation, which we omit.
\end{proof}

Since it plays a relevant role in what follows, we state a result of Exel \cite{E} for the reader's convenience.

%
\begin{theorem}[{\cite[Theorem V.11]{E}}]\label{Exeltheo}
Let $A$ be a unital C*-algebra, let $\alpha$ be an automorphism of $A$  and let $\tau$ be an $\alpha$-invariant trace on $A$.  If  $n\geq 1$ and $u$ is a smooth path 
in $U_n(M\sb\alpha)=\{ f\colon [0,1]\rightarrow U_n(A) \; |\; f(1)=\alpha(f(0))\}$ then
\[
   \tau( k\sb\alpha( [u] )) = \frac{1}{2\pi i} \int_{0}^{1} \tau \left( u\sp\prime (s)\sp\ast u(s)\right)
                               \, \d\lambda (s)
\]
where $k_{\alpha}\colon K_1(M\sb\alpha)\rightarrow K_0(\cprod{A}{\alpha})$ is the Connes' Isomorphism (\ref{connesk}) and $\lambda$ is  Lebesgue measure on $[0,1]$.
\end{theorem}

The following proposition is the key to the main result of this section.

\begin{proposition}\label{ttauvstau}
Let the notation  be as in Definition~\ref{ttau}. Let  $t\in\R$ and let  $M_{\phi_t}$ be the mapping torus of  $\phi_t$. Then the following diagram is commutative.
\[
   \xymatrix{ %
             K_1(M\sb\alpha) \ar[rr]^{\theta_t} \ar[d]_{k_{\alpha}}  & & K_1(M_{\phi_t}) \ar[d]^{k_{\phi_t}} \\
             K_0(\cprod{A}{\alpha}) \ar[dr]_{t\tau}  &  & K_0(\cprod{M\sb\alpha}{\phi_t}) \ar[dl]^{\ttau}  \\
                                                     & \R  &    \\
                  }%
\]
\end{proposition}
\begin{proof}
Let $x\in K_1(M\sb\alpha)$. Choose $u$ to be a smooth path in $U\sb\infty (M\sb\alpha)=\{ f\colon [0,1]  \rightarrow U\sb\infty (A)\; |\; f(1)=\alpha( f(0))\}$ 
such that  $u$ is constant on small neighborhoods of $0$ and $1$ and $x=[u]$ (such a choice is always possible since any element in $U\sb\infty(M\sb\alpha)$ is homotopic to one of these).
Then
\begin{eqnarray*}
      ( \ttau \circ k_{\phi_t}\circ\theta_t ))(x)  & = & 
                           (\ttau \circ k_{ \phi_t } )( [ \phi_{ st }( u ) ]  )
                             \,\,\text{ by Definition~\ref{theta}} \\
                 & = & \frac{1}{2\pi i} \int_{0}^{1} \ttau \left(  \left({\textstyle{\frac{\d}{\d s}}}\phi_{st} (u) \right) \sp\ast \phi_{st}(u)\right) \,\d\lambda( s)
                              \,\,\text{ by Theorem~\ref{Exeltheo}} \\
                 & = & \frac{1}{2\pi i} \int_{0}^{1} t\,\ttau \left( \phi_{st}\left(\left(u\sp\prime\right)\sp\ast\right) \phi_{st} (u) \right) \,\d\lambda( s )
                              \,\,\text{ by Lemma~\ref{derivativeofphi}}\\
                 & = & \frac{t}{2\pi i} \int_{0}^{1} \ttau \left(\left(u\sp\prime\right)\sp\ast u \right)\, \d\lambda( s )
                              \,\,\text{ by $\phi$ invariance of $\ttau$}\\
                 & = & \frac{t}{2\pi i}\,\ttau\left(  \left(u\sp\prime\right)\sp\ast  u \right)
                              \,\,\text{by computation of the integral}\\
                  & = & \frac{t}{2\pi i}\int_{0}^{1} \tau \left(u\sp\prime (t)\sp\ast u(t) \right)\, \d\lambda(t) 
                               \,\,\text{ by Definition~\ref{ttau}} \\
                  & = & (t\,\tau\circ k_{\alpha} )\left( x \right)
                               \,\,\text{ by Theorem~\ref{Exeltheo}},
   \end{eqnarray*}
as was to be proved.
\end{proof}

%
%

We are ready to prove the main result of this section.
\begin{theorem}\label{mainsection1}
Let $A$ be a unital C*-algebra, let $\alpha$ be an automorphism of $A$ and let $\tau$ be an $\alpha$-invariant trace on $A$. Consider the mapping torus flow $(M\sb\alpha,{\phi})$ associated to $(A,\alpha)$ and let $\ttau$ be the  ${\phi}$-invariant trace on $M\sb\alpha$ induced by $\tau$, as in Definition~\ref{ttau}. If $t\in\R$ then
\[
  K_0(\cprod{M\sb\alpha}{\phi_t})\cong K_0(M\sb\alpha)\oplus K_0(\cprod{A}{\alpha})\cong
  K_1(\cprod{M\sb\alpha}{\phi_t})
\]
and, for each $(x,y)\in K_0(M\sb\alpha)\oplus K_0(\cprod{A}{\alpha})$,
\[
   \ttau(x,y) = \ttau (x) + t\,\tau (y).
\]
\end{theorem}
\begin{proof}
Since $\phi_t$ is homotopic to the identity map $\id$ on $M\sb\alpha$, it follows that the $K$-groups of $\cprod{M\sb\alpha}{\phi_t}$ and $\cprod{M\sb\alpha}{\id}$ are isomorphic (see, e.g.\ \cite[Corollary~10.5.2]{Black}). But since the $K$-groups of $\cprod{M\sb\alpha}{\id}$  are isomorphic to $K_0(M\sb\alpha)\oplus K_1(M\sb\alpha)$ (cf.\ \cite[10.1.1(a) and 9.4.1]{Black}) and
$K_1(M\sb\alpha)$ is isomorphic to $K_0(\cprod{A}{\alpha})$ (via Connes' Isomorphism (\ref{connesk})),
we conclude that 
\[
 K_0(\cprod{M\sb\alpha}{\phi_t})\cong K_0(M\sb\alpha)\oplus K_0(\cprod{A}{\alpha})\cong
  K_1(\cprod{M\sb\alpha}{\phi_t}),
\]
as desired.
Now, from Proposition~\ref{ttauvstau} we obtain the following commutative diagram.
\[
   \xymatrix{ %
  K_0(\cprod{A}{\alpha})\ar[rr]^(.49){k_{\phi_t}\circ\theta_t\circ k_{\alpha}^{-1} } \ar[drr]_{t\,\tau} 
                                            &&  K_0(\cprod{M\sb\alpha}{\phi_t})\ar[d]^{\ttau} 
                                            &&  K_0(M\sb\alpha)\ar[ll]_(.4){i\sb\ast} \ar[dll]^{\ttau} \\
                                             && \R &  & \\
                 }
\]

An application of  Theorem~\ref{kgroup} completes the proof.
\end{proof}

\section{Entropy}
Topological entropy has been very successful in topological dynamics as an invariant for conjugacy of dynamical systems. Roughly speaking, it measures the total exponential complexity of the orbit structure with a single number 
(see, e.g.\ \cite[Chapter~7]{W}). Voiculescu introduced a notion of entropy for automorphisms of unital nuclear C*-algebras based on local approximation \cite{Vo} and Brown subsequently extended this notion to automorphisms of exact C*-algebras using nuclear embeddability {\cite{Nate}}. The Voiculescu-Brown entropy is an extension of the original notion of entropy in topological dynamics. Indeed, it was proved in \cite[Proposition~4.8]{Vo}  that the topological entropy for a homeomorphism on a compact metric space coincides with the Voiculescu-Brown entropy for the automorphism (of the commutative C*-algebra) induced by the homeomorphism.

Suppose that $A$ is an exact C*-algebra and $\alpha$ is an automorphism of $A$.
Denote by $(M\sb\alpha,{\phi})$ the mapping torus flow associated to  $(A,\alpha)$.  Then $M\sb\alpha$ is an exact C*-algebra and $\phi_t$ is an automorphims of $M\sb\alpha$ for all $t\in\R$. ($M\sb\alpha$ is exact because it is (isomorphic to) a subalgebra of $C([0,1])\otimes A$ which is exact since it is the tensor product of a nuclear and an exact C*-algebra.) In this section we find a formula to compute the Voiculescu-Brown entropy of $\phi_t$ in terms of $t$ and the Voiculescu-Brown entropy of $\alpha$.

We start by recalling the definition of the Voiculescu-Brown entropy \cite{Nate}.
Let $A$ be an exact C*-algebra and let $\alpha$ be an automorphism of $A$. Let 
$\pi\colon A\rightarrow \calBH$ be a faithful *-representation of $A$.  We denote by $CPA(\pi,A)$ the collection of all triples $(\phi,\psi,B)$ where $B$ is a finite dimensional C*-algebra and $\phi\colon A\rightarrow B$ and 
$\psi\colon B\rightarrow \calBH$ are contractive completely positive maps.
 Given $\delta>0$ and a finite set $\omega\subset A$,  we denote by 
 \begin{equation}\notag
 \begin{split}
 rcp(\omega,\delta)=\inf\{\mathrm{rank}(B) \colon (\phi,\psi,B)\in CPA(\pi,A)\text{ and }  \\ 
   \| \psi\circ\phi (x) - \pi(x)  \|<\delta  
   \text{ for all } x\in\omega \}
 \end{split} 
\end{equation}
where $rank(B)$ is the dimension of a maximal abelian subalgebra of $B$. As the notation indicates, $rcp(\omega,\delta)$ is independent of the faithful *-representation $\pi$, as proved in \cite[Proposition~1.3]{Nate}. 

We now define
\begin{itemize}
  \item $\entr (\alpha,\omega,\delta)={\ds{\limsup_{m \rightarrow\infty}}}\,
                \frac{1}{m} \log( rcp (\omega\cup \alpha(\omega)\cup\cdots\cup \alpha^{m-1}(\omega),
                \delta))$
   \item $\entr (\alpha,\omega)={\ds{\sup_{\delta>0}}}\, \entr(\alpha,\omega,\delta) $
   \item $\entr( \alpha )= {\ds{\sup_{\omega}}}\, \entr(\alpha,\omega)$
\end{itemize}
where the last supremum is taken over all finite sets $\omega\subset A$.  We call this last quantity the {\em Voiculescu-Brown entropy} of $\alpha$.

The following proposition extends a result of Bowen \cite[Proposition~21]{Bowen} to the non-commutative setting. Recall that an action of $\R$ on a C*-algebra $A$ is a homomorphism $\phi$ of $\R$ into the group of automorphisms of $A$ such that $s\mapsto \phi_s(a)$ is continuous for each fixed $a\in A$. 

\begin{proposition}\label{entr1}
Let $A$ be an exact C*-algebra and let ${\phi}$ be an action of $\R$ on $A$.  If $t\in\R$ then
\[
    \entr( \phi_t) =|t|\,\entr (\phi_1).
\]
\end{proposition}
\begin{proof}
To prove the proposition, it will suffice to show that
\begin{equation}\label{bowen}
 \entr(\phi_t)\leq {\ds{\frac{t}{s}}}\, \entr(\phi_s)
\end{equation}
for all reals $s,t>0$. Indeed, this would give us that the proposition is valid for all reals $t>0$. This would then imply the proposition for reals $t<0$ since by \cite[Proposition~2.5]{Nate} we also get
\[
\entr(\phi_t)=\entr(\phi_{-t}^{-1})=\entr(\phi_{-t})=-t\,\entr(\phi_1)=|t|\, \entr(\phi_1).
\]

We now prove (\ref{bowen}). For this purpose, let $s,t>0$ be real numbers. 
Let $\omega\subset A$ be a finite set and let $\epsilon >0$.
Since the set
\[
\Gamma=\{ \phi_r(a)\colon 0\leq r\leq s,\,\, a\in\omega\}
\] 
is compact, it can be covered by finitely many $\ts{\frac{\epsilon}{4}}$-balls with centers in $\Gamma$. Let 
{\allowbreak $\phi_{r_1}(a_1),\phi_{r_2}(a_2),$}{ $\ldots,\phi_{r_k}(a_k)$} be the centers of such balls.
Put
\[
\omega\sp\prime=\{\phi_{r_i}(a_i)\colon 1\leq i\leq k\}.
\]

By \cite[Proposition~1.3]{Nate}, we may assume that $A$ is faithfully represented on a Hilbert space $\calH$. Let  $(\alpha,\beta, B)\in CPA (\id_A ,A)$.
Then  $\|\beta\circ\alpha (b) - b\| <\epsilon$ for all $b\in\Gamma$ whenever
$\| \beta\circ\alpha (a) -  a\|<\ts{\frac{\epsilon}{2}}$ for all $a\in\omega\sp\prime$.
Let $\delta=\ts{\frac{\epsilon}{2}}$. This proves that 
\[
   rcp(\omega,\epsilon)\leq rcp(\omega\sp\prime,\delta).
\]
Furthermore, since $\phi_s$ is an isometry  we obtain, whenever $(m-1)t<
(n-1)s$, that 
\[
  rcp(\omega\cup\phi_t(\omega)\cup\cdots\cup\phi_{t}^{m-1}(\omega),\epsilon)
  \leq
  rcp(\omega\sp\prime\cup\phi_s(\omega\sp\prime)\cup\cdots\cup
  \phi_{s}^{n-2}(\omega\sp\prime),\delta).
\]
Hence
\begin{eqnarray*}
   \entr(\phi_t,\omega,\epsilon)&\leq&
       \ds{\limsup_{m\rightarrow\infty} }\,
       {\ts{\frac{1}{m}}} \,\log\left( rcp  
       \left(\omega\sp\prime\cup\phi_s(\omega\sp\prime)\cup\cdots\cup
  \phi_{s}^{\ts{\left[ \ts{{(m-1)t}/{s}}\right] }}(\omega\sp\prime),\delta\right)\right)\\
      &\leq &
      \entr(\phi_s,\omega\sp\prime,\delta)\, \ds{\limsup_{m\rightarrow\infty}}\,
      \ts{\frac{1}{m}}  \left( \left[ (m-1)t/s \right]  \right) \\
      &=&
       {\frac{t}{s}}\,\entr(\phi_s,\omega\sp\prime,\delta)\\
      &\leq&
      {\frac{t}{s}} \,\entr(\phi_s). 
\end{eqnarray*}

Letting $\epsilon$ and $\omega$ vary, $\entr(\phi_t)\leq \ts{\frac{t}{s}}\entr(\phi_s)$,
as wanted.
\end{proof}

The following proposition is a slight generalization of \cite[Proposition~2.10]{Nate}, that is, without assuming that the C*-algebras involved  are unital.  As N.\ Brown kindly pointed out to me, a proof can be obtained by imitating the proof of \cite[Proposition~2.10]{Nate} without making any changes.

\begin{proposition}[N.\ Brown]\label{Nate}
  Let $A$ and $B$ be two exact C*-algebras and let  $\alpha$ and $\beta$ be automorphisms of $A$ and $B$, respectively. If $\rho\colon A\to B$ is an injective completely positive map such that $\rho^{-1}\colon \rho(A)\to A$ is also completely positive and $\rho\circ \alpha = \beta\circ\rho$, then $\entr(\alpha)\leq \entr(\beta)$.
\end{proposition}

Let $A$ be an exact C*-algebra and let  $\alpha$ be an automorphism of $A$. Let $(M\sb\alpha,\phi)$ be the mapping torus flow associated to $(A,\alpha)$. Then, for $f\in M\sb\alpha$, $\phi_1(f)=\alpha\circ f$.  It is then expected that $\phi_1$ and $\alpha$ have the same Voiculescu-Brown entropy. We prove this in the following.

\begin{proposition}\label{entr2}
Let $A$ be an exact C*-algebra and let $\alpha$ be an automorphism of $A$.
Consider the mapping torus flow $(M\sb\alpha,{\phi})$ associated to $(A,\alpha)$.
Then
\[
   \entr( \phi_1) =  \entr(\alpha).
\]
\end{proposition}
\begin{proof}
For each $a\in A$, define $f_a\in M\sb\alpha$ by the formula $f_a(s)=(1-s)a+s\alpha(a)$, for $s\in[0,1]$. Let $\rho\colon A\to M\sb\alpha$ be the linear map given by $\rho(a)=f_a$. Using the fact that $\alpha$ is an automorphism and that the set of positive elements of a C*-algebra forms a closed cone, we see that $\rho$ is a completely positive map. We claim that also $\rho^{-1}\colon \rho(A)\to A$ is completely positive. To prove that $\rho^{-1}\colon \rho(A)\to A$ is positive, one considers a positive element $f_a$ in $\rho(A)$. Then $f_a=g\sp\ast g$ for some $g\in M\sb\alpha$. Hence 
\[g\sp\ast(0)g(0)=f_a(0)=a=\rho^{-1}(f_a)
\]
is positive, as wanted. The proof of the claim will be completed after we consider the isomorphism between  $M\sb\alpha(A)\otimes M_{n}$ and
$M\sb{\alpha\otimes \id_n}(A\otimes M_n)$ and repeat the argument above;  we omit  details. Futhermore, since $\rho\circ\alpha=\phi_1\circ\rho$, we apply Proposition~\ref{Nate} to obtain
$\entr(\alpha)\leq\entr(\phi_1)$.

For opposite inequality, consider the inclusion map $\iota\colon M\sb\alpha\rightarrow C([0,1],A)$. Then it is clear that both $\iota$ and $\iota^{-1}\colon \rho(M\sb\alpha)\to M\sb\alpha$ are completely positive.  Consider now $\psi\in \Aut(C([0,1],A))$ defined by $\psi(f)=\alpha\circ f$ for $f\in C([0,1],A)$. Then $\iota$ intertwines $\phi_1$ and $\psi$, that is, $\iota\circ\phi_1=\iota\circ\psi$. Using Proposition~\ref{Nate} we conclude that 
\begin{equation}\label{1}
\entr(\phi_1)\leq \entr(\psi).
\end{equation}
Now, since $\psi$ is conjugate to $\id\otimes\alpha \in\Aut\left(C([0,1])\otimes A\right)$, it follows that
\begin{equation}\label{2}
  \entr(\psi)=\entr(\id\otimes\alpha).
\end{equation}
But since $\entr(\id)=0$, by \cite[Proposition~2.7]{Nate}, we also get
\begin{equation}\label{3}
   \entr(\id\otimes\alpha)=\entr(\alpha).
\end{equation}
Thus, combining (\ref{1}), (\ref{2}) and (\ref{3}) we obtain $\entr(\phi_1)\leq\entr(\alpha)$. This completes the proof.
\end{proof}

We are ready to prove the main result of this section.
\begin{theorem}\label{mainsection2}
Let $A$ be an exact C*-algebra and let $\alpha$ be an automorphism of $A$.
Suppose that $(M\sb\alpha,{\phi})$ is the mapping torus flow corresponding to $(A,\alpha)$.
If $t\in\R$ then
\[
   \entr (\phi_t) =|t|\, \entr(\alpha).
\]
\end{theorem}
\begin{proof}
   The proof of this theorem is just an application of Propositions \ref{entr1} and \ref{entr2}.
\end{proof}


\section{Applications }

  In this section we specialize to the commutative setting  some of the  results of the previous sections. We prove that crossed products by minimal time-$t$ homeomorphisms of suspensions built over strongly orbit equivalent Cantor minimal systems have isomorphic Elliott invariants. As an application, we show the existence of dynamical systems on (compact metric) connected 1-dimensional spaces which are not flip conjugate (because of different entropy) yet their associated crossed products have isomorphic Elliott invariants. We will also address some interesting questions emanating from our investigation.

  We begin with some definitions complementing the ones given in the introduction.  Let $(X,S)$ be a dynamical system. We say that $(X,S)$ (or just $S$) is {\em minimal} if there is no non-trivial $S$-invariant closed subset of $X$.  Equivalently, $(X,S)$ is minimal if the orbit $\{S^n(x)\colon n\in\Z\}$ of each $x\in X$ is dense in $X$. By $M(X,S)$ we will denote the set of $S$-invariant Borel probability measures. If $M(X,S)$ has only one element, then we say that $(X,S)$ (or just $S$) is {\em uniquely ergodic}.  When a dynamical system is both minimal and uniquely ergodic, we say that it is {\em strictly ergodic}. If $T$ is an action of $\R$ on the compact metric space $Y$, we denote by $M(Y,T)$ the set of all $T$-invariant Borel probability measures on $Y$  (i.e.\ $M(Y,T)$ contains all Borel probability measures $\nu$ on $Y$ satisfying $\nu(T^t(E))=\nu(E)$ for all $t\in\R$ and for all measurable sets $E\subset Y$).
In analogy  to our terminology for (discrete) dynamical systems, we say that $(Y,T)$ is minimal, uniquely ergodic or strictly ergodic  if there is no non-trivial closed $T$-invariant subset of $Y$, $M(Y,T)$ has only one element or is both minimal and uniquely ergodic, respectively. 
  
   It is well known that if $(X,S)$ is a minimal (uniquely ergodic)  dynamical system and $X$ is infinite then its associated crossed product $C\sp\ast(\Z,X,S)$ is simple (has a unique normilized trace). For this reason, we will always assume $X$ to be infinite.
  
  Two dynamical systems $(X_1,S_1)$ and $(X_2,S_2)$ are said to be 
  \begin{itemize}
     \item {\em conjugate} if there is a homeomorphism $F\colon X_1\rightarrow X_2$, called a conjugation map, such that $F\circ S_1=S_2\circ F$.  
     \item {\em flip conjugate} if $(X_1,S_1)$ is conjugate to either $(X_2,S_2)$ or $(X_2,S_{2}^{-1})$.
     \item {\em orbit equivalent} if there is  homeomorphism $F\colon X_1\rightarrow X_2$, called an orbit map, such that $F(\{ S_{1}^{n}(x)\colon n\in\Z\}) = \{ S_{2}^{n}(F(x))\colon
n\in\Z\}$ for all $x\in X_1$. If in addition $(X_1,S_1)$ and $(X_2,S_2)$ are minimal then we can uniquely define two functions $m,n\colon X_1\rightarrow \Z$, called the orbit cocycles of $F$, such that $F(S_1(x))=S_{2}^{m(x)}(F(x))$ and $F(S_{1}^{n(z)}(x))=
S_2(F(x))$ for $x\in X_1$. When there is an orbit map so that the associated orbit cocycles have at most one point of discontinuity each, then we say that $(X_1,S_1)$ and $(X_2,S_2)$ are {\em strongly orbit equivalent}, cf.\ \cite[Definition~1.3]{GPS}.
  \end{itemize}
  
  We now proceed to describe the Elliott invariants for crossed products associated to minimal dynamical systems.     Let $A$ be the crossed product associated to a minimal dynamical system $(X,S)$.  Let $T(A)$ denote the set of (normalized) traces on $A$, let $S(K_0(A))$ be the state space of the group $K_0(A)$, and let $r_A\colon T(A)\rightarrow S(K_0(A))$ be the canonical map given by $r_A(\tau)(w)=\tau(w)$ for $\tau\in T(A)$ and $w\in K_0(A)$. The Elliott invariant of $A=C\sp\ast(\Z,X,S)$ is given by the following data. 
\begin{itemize}
    \item The abelian group $K_1(A)$, which we will denote by $K^1(X,S)$.
    \item The scaled ordered abelian group $K_0(A)$, which we will denote by $K^0(X,S)$.  The scale is the distinguished element $[1]$ and the order is defined by $x\geq0$ if and only if $x=[p]$ for some integer $n\geq 1$ and some projection $p\in M_n(A)$.
    \item The simplex $T(A)$ equipped with the weak* topology, cf.\ \cite[Theorem~3.1.18]{Sakai}.
    \item The natural continuous affine map $r_A\colon T(A)\rightarrow S(K_0(A))$.
\end{itemize} 
  The Elliott invariant of $A$ is, according to the {\em Elliott conjecture}, a complete isomorphism invariant for $A$.
   
 An {\em isomorphism} from the Elliott invariant of  $A_1=C\sp\ast(\Z,X_1,S_1)$ to that of  $A_2=C\sp\ast(\Z,X_2,S_2)$ consists of an isomorphism 
 $\phi_1\colon K^1(X_1,S_1)\rightarrow K^1(X_2,S_2)$ of (abstract) abelian groups,
  an isomorphism $\phi_0\colon K^0(X_1,S_1)\rightarrow K^0(X_2,S_2)$ of scaled ordered abelian groups, and an affine homeomorphism $f\colon T(A_2)\rightarrow T(A_1)$ 
  such that the diagram
\[  
  \xymatrix{  %
     T(A_2)          \ar[r]^f \ar[d]_{r_{A_2}}  & T(A_1) \ar[d]^{r_{A_1}}    \\
    S(K_0(A_2))  \ar[r]_{\widehat{\phi}_0}                        &S(K_0(A_1))   \\                                             
                }
  \] 
  commutes, where $\widehat{\phi}_0$ is the dual of $\phi_0$ (i.e. $\widehat{\phi}_0(\sigma)=\sigma\circ\phi_0$ for $\sigma\in S(K_0(A_2))$).
   In this case we will say that $A_1$ and $A_2$ have {\em isomorphic Elliott invariants.}

\medskip
We will need the following natural result, which is most  likely well known to the specialist. We could not find the second part of it in the literature, however. We sketch here a proof for completeness.

 \begin{lemma}\label{originaliffsuspension}
 Let $(X,S)$ be a dynamical system and let $(Y,T)$ be its associated suspension flow.
 Then we have the following.
 \begin{enumerate}
   \item $(X,S)$ is minimal if and only if $(Y,T)$ is minimal.
   \item There is an affine homeomorphism between $M(X,S)$ and $M(Y,T)$.
   \item $(X,S)$ is strictly ergodic if and only if $(Y,T)$ is strictly ergodic. 
 \end{enumerate}
 \end{lemma}
 
 \begin{proof}
    Part (1)  is \cite[5.12.12]{V}. Part (3) follows trivially from part (1) and (2). We now sketch the proof of part (2). 
    
   If $\mu$ is in $M(X,S)$ then $\mu$ induces an element $\nu$ in $M(Y,T)$ via the formula
   \[
         E \mapsto  (\mu\times \lambda)\left(\pi^{-1} (E) \right)
   \]
where $\lambda$ is Lebesgue measure on $[0,1]$, $E$ is a measurable set in $Y$, and $\pi$ is the canonical quotient map $X\times [0,1]\rightarrow Y$.
Conversely, suppose that $\nu$ is a measure in $M(Y,T)$. Define $\mu$ by
\begin{equation*}
         F\mapsto {\nu}(\pi(F\times [0,1))
\end{equation*}
where $F$ is a measurable set in $X$.  One checks that  $\mu$ is in $M(X,S)$ 
and that this correspondence defines an affine homeomorphism between $M(X,S)$ and $M(Y,T)$, as desired.
 \end{proof}

   We now combine Theorem~\ref{mainsection1} with \cite[Theorem~4.5]{NCP2} to prove a result which will allow us to compute the Elliott invariants of crossed products by minimal time-$t$ homeomorphisms on suspensions built over  finite dimensional dynamical systems.

\begin{proposition}\label{elliottinvariants}
Let $(X,S)$ be a dynamical system and let $(Y,T)$ be the suspension flow  associated to 
$(X,S)$. 
 Suppose that $\tau\sb\mu$ and $\tau_{\nu}$ are the traces on $C(X)$ and $C(Y)$ associated to some $\mu\in M(X,S)$ and its corresponding $\nu\in M(Y,T)$, respectively. If $t\in \R$  then 
\[
K^1(Y,T^t)\cong K^0(Y)  \oplus K^0(X,S)\cong K^0(Y,T^t)
\]
  and, for each $(w,z)\in K^0(Y)\oplus K^0(X,S)$,
\[  
        \tau_{\nu}(w,z )=\tau_{\nu}(w) + t\, \tau\sb\mu(z).
\]
Furthermore, if $X$ is finite dimensional and $T^t$ is minimal then $A=C\sp\ast(\Z,Y,T^t)$ satisfies the following K-theoretical version of Blackadar's Second Fundamental Comparability Question: if $w\in K^0(Y,T^t)$ satisfies $\tau(w)>0$ for all (normalized) traces $\tau$ on $A$, then there is a projection $p\in M\sb\infty (A)$ such that $x=[p]$.
\end{proposition}
\begin{proof}
    Let $(M\sb\alpha,{\phi})$ be the mapping torus flow associated to $(C(X),\alpha)$, where $\alpha$ is the automorphism of $C(X)$ induced by $S$. It follows that the automorphism of $C(Y)$ induced by  $T^t$ is conjugate to $\phi_t$ using the conjugation map  $\sigma \colon M\sb\alpha \rightarrow C(Y)$ given by $\sigma (f)\left([x,s]\right)=(f(s))(x)$,
where  $f\in M\sb\alpha$ and $[x,s]\in Y$. Therefore
 $\sigma$ induces an isomorphism  between $\cprod{M\sb\alpha}{\phi_t}$ and $C\sp\ast (\Z,Y,T^t)$. Let $\ttau$ be the trace on $M\sb\alpha$ induced by $\tau\sb\mu$ as in Definition~\ref{ttau}. One also verifies  that  $\tau_{\nu}\circ\sigma=\ttau$. Hence, the proposition follows from Theorem~\ref{mainsection1} by taking $A$, $\alpha$ and $\tau$ to be $C(X)$, the automorphism induced by $S$ and $\tau\sb\mu$, respectively.

If $X$ is finite dimensional and $T^t$ is minimal then $Y$ is also finite dimensional and so we may apply \cite[Theorem~4.5]{NCP2} to obtain the last assertion of the corollary.
\end{proof}

   It has been proved by Giordano, Putnam and Skau \cite[Theorem~2.1]{GPS} that strong orbit equivalence of Cantor minimal systems corresponds to isomorphism  of their associated crossed products. In the following theorem we show that  crossed products by minimal  time-$t$ homeomorphisms of suspensions built over strongly orbit equivalent Cantor minimal systems have isomorphic Elliott invariants.
   
 \begin{theorem}\label{ElliottInv}
 Let $(X_1,S_1)$ and $(X_2,S_2)$ be strongly orbit equivalent Cantor minimal systems and let  $(Y_1,T_1)$ and $(Y_2,T_2)$ be the suspension flows associated to $(X_1,S_1)$ and 
 $(X_2,S_2)$, respectively.  Let $t\in\R$ such that
 \begin{enumerate}
    \item[(i)] $T_{1}^{t}$ and $T_{2}^{t}$ are minimal,
    \item[(ii)]  $M(Y_1,T_1)=M(Y_1,T_{1}^{t})$ and $M(Y_2,T_2)=M(Y_2,T_{2}^{t})$. 
\end{enumerate}
 Then  $C\sp\ast(\Z,Y_1,T_{1}^{t})$ and $C\sp\ast(\Z,Y_2,T_{2}^{t})$ have isomorphic Elliott invariants.
 \end{theorem}  
\begin{proof}
  Since $(X_1,S_1)$ and $(X_2,S_2)$ are strongly orbit equivalent, one may assume that $S_1$ and $S_2$ act on the same space $X$, have the same orbits and that the orbit cocycles $m,n\colon X\rightarrow \Z$ defined by $S_1(x)=S_{2}^{n(x)}(x)$ and
  $S_2(x)=S_{1}^{m(x)}(x)$, for all $x\in X$, have at most one point of discontinuity. It also follows that $M(X,S_1)=M(X,S_2)$, cf.\  the proof of ($i$)$\Rightarrow$($ii$) in \cite[Theorem~ 2.1]{GPS}, pages 78--79.
  For $i=1,2$, we may identify $K^0(X,S_i)$ with $C(X,\Z)\slash \im\left(\id-(S_i)\sb\ast\right)$ where $C(X,
\Z)$ denotes the (countable) abelian group (under addition) of continuous functions on $X$ with values in $\Z$ and $\im\left(\id-(S_i)\sb\ast\right)=\{f-f\circ S_{i}^{-1}\colon f\in C(X,\Z\}$ (see \cite[Theorem~1.1]{Putnam}).  Using this identification, let $[f]_i$, where $f\in C(X,S_i)$, denote the class in $K^0(X,S_i)$ represented by $f$. Following the proof  ($i$)$\Rightarrow$($ii$) in \cite[Theorem~2.1]{GPS}, pages 78--79, we see that the map $I\colon K^0(X,S_1)\rightarrow K^0(X,S_2)$ defined by $I([f]_1)=[f]_2$ for $f\in C(X,\Z)$, gives an scaled ordered isomorphism of abelian groups.

On the other hand, for $i=1,2$, since $K^1(X,S_i)$ is isomorphic to $\Z$ (cf.\ \cite[Theorem~1.1]{Putnam}) and since $C(Y_i)$ is isomorphic to the mapping torus of the automorphism of $C(X)$ induced by $S_i$, we may use the Connes' Isomorphism (\ref{connesk}) to obtain that $K^0(Y_i)$ is isomorphic to $\Z$ and so $K^0(Y_i)$ is generated by $[ 1]$. Thus, from Proposition~\ref{elliottinvariants}, we conclude that $K^1(Y_1,T_{1}^{t})$ is isomorphic to $K^1(Y_2,T_{2}^{t})$.

Let us denote by  $F_i\colon M(X,S_i)\to M(Y_i,T_i)$, for $i=1,2$,
the affine homeomorphisms obtained by Lemma~\ref{originaliffsuspension}.
Define a function $\mathcal{F}\colon T(C\sp\ast (\Z,Y_2,T_{2}^{t})) \rightarrow T(C\sp\ast (\Z, Y_1,T_{1}^{t}))$ in the following natural way. 
Let  $\tau_{\nu_2}$ be in $T(C\sp\ast (\Z,Y_2,T_{2}^{t}))$ where 
$\nu_2$ is an element of $M(Y_2,T_{2}^{t})$.
 Then since $M(X,S_1)=M(X,S_2)$, $M(Y_1,T_1)=M(Y_1,T_{1}^{t})$ and $M(Y_2,T_2)=M(Y_2,T_{2}^{t})$, 
we may let
$\mathcal{F}(\tau_{\nu_2})$  be $\tau_{\nu_1}$, where $\nu_1=(F_1\circ F_{2}^{-1})(\nu_2)$. Since $F_1$ and $F_2$ are affine homeomorphisms,  one verifies that $\mathcal{F}$ is also an affine homeomorphism.

We claim that $\id\oplus I\colon \Z\oplus K^0(X,S_1)\rightarrow \Z\oplus K^0(X,S_2)$ gives a scaled ordered isomorphism from $K^0(Y_1,T_{1}^{t})$ onto $K^0(Y_2,T_{2}^{t})$. Indeed, suppose that $0<(n,[f]_1)\in\Z\oplus K_0(X,S_1)$. Let  $\nu_2$ be a measure in $M(Y_2,T_{2}^{t})$ and let $\mu=F_{2}^{-1}(\nu_2)$ and $\nu_1=F_1(\mu)$ be the corresponding measures in $M(X,S_1)$ and $M(Y_1,T_{1}^{t})$. Then, using Proposition~\ref{elliottinvariants} in the second and third step and the fact that $(n,[f]_1)>0$ and $C\sp\ast(\Z,Y_1,T_{1}^{t})$ is simple in the last step below, we get:
\begin{eqnarray*}
    \tau_{\nu_2} \left(\left(\id \oplus I\right)\left(n,[f]_1\right)\right) &=& \tau_{\nu_2} \left(n,[f]_2\right) \\
             &=&  n+ t\tau\sb\mu \left( f \right)  \\
             &=& \tau_{\nu_1}(n,[f]_1)\\
             & > & 0.
\end{eqnarray*}
Since $\nu_2$ was arbitrary, the last part of Proposition~\ref{elliottinvariants} gives us that $(\id\oplus I)(n,[f]_1)>0$. Hence $\id\oplus I$ is an scaled ordered isomorphism, proving the claim.

To complete the proof of the theorem, let $\tau_{\nu_2}$ be in $M(Y_2,T_{2}^{t})$ where $\nu_2$ is an element of $M(Y_2,T_{2}^{t})$ and let
$(n,[f]_1)$ be in $\Z\oplus K^0(X,S_1)$. If  $\mu=F_{2}^{-1}(\nu_2)$ and $\nu_1=F_1(\mu)$ are the measures in $M(X,S_1)$ and $M(Y_1,T_{1}^{t})$ corresponding to $\nu_2$ then, by using Proposition~\ref{elliottinvariants} in the second and third step and the definition of $\mathcal{F}$ in the last step below, we get:
\begin{eqnarray*}
     \left(\widehat{\id\oplus I}\circ r_{A_2}\right) (\tau_{\nu_2})( n,[f]_1) &=& r_{A_2}(\tau_{\nu_2})(n,[f]_2)\\
          &=& n +t\tau\sb\mu (f)\\
          &=& r_{A_1}(\tau_{\nu_1})(n,[f]_1)\\
          &=& (r_{A_1}\circ \mathcal{F})(\tau_{\nu_2})(n,[f]_1).
\end{eqnarray*}
  Hence  $C\sp\ast(\Z,Y_1,T_{1}^{t})$ and $C\sp\ast(\Z,Y_2,T_{2}^{t})$ have isomorphic Elliott invariants.
\end{proof}

We now turn to the problem of finding minimal time-$t$ homeomorphisms of suspensions.  Let $(X,S)$ be a dynamical system and let $(Y,T)$ be the suspension flow associated to $(X,S)$.
The following proposition will tell us that the minimality of $(X,S)$ is a necessary and sufficient condition for the existence of $t\in\R$ such that $T^t$ is minimal.  The author would like to thank E.\ Glasner for referring us to a result in his book \cite{G} which helped to improve our original proposition. 
  
 \begin{proposition}\label{texists}
 Let $(X,S)$ be a dynamical system and let $(Y,T)$ be the suspension flow associated to $(X,S)$. Then 
 \begin{enumerate}
 \item$(X,S)$ is minimal if and only if the set
         \[ \left\{t\in \R\colon T^t\text{ is {\rm not} minimal}\right\} \]
         is a countable subset of $\R$ which contains $\Q$.
 \item  $(X,S)$ is strictly ergodic if and only if the set
         \[ \left\{t\in \R\colon T^t\text{ is {\rm not} strictly ergodic}\right\} \]
         is a countable subset of $\R$ which contains $\Q$.
 \end{enumerate} 
 \end{proposition}
 \begin{proof}
 If there is  $t\in\R$ such that $T^t$ is minimal  then $S$ must be minimal. (If $M\subset X$ is closed nontrivial $S$-invariant then $\pi (M\times [0,1])\subset Y$ is closed nontrivial $T^t$-invariant.)  For the converse, observe that $T^r$ is not minimal whenever $r$ is rational.  Assume now that $(X,S)$ is minimal. 
 Then $(Y,T)$ is minimal by Lemma~\ref{originaliffsuspension} and therefore
  \[
     \{t\in \R\colon T^t\text{ is {\rm not} minimal}\}
 \]
   is a countable set, cf. \cite[4.24.1]{G}.

For the second part of the proposition, suppose that $T^t$ is strictly ergodic. Then $S$ must also be strictly ergodic. (An argument for minimality is given in the above paragraph. For unique ergodicity, observe that if $M(X,S)$ has more that one element then so does $M(Y,T)$ by Lemma~\ref{originaliffsuspension}. Since $M(Y,T)\subset M(Y,T^t)$ then $M(Y,T^t)$ also has more than one element.) For the converse, observe that $T^r$ is not strictly ergodic whenever $r$ is rational. Assume now that $(X,S)$ is strictly ergodic.
Then $(Y,T)$ is strictly ergodic by  Lemma~\ref{originaliffsuspension} and therefore
 \[
     \{t\in \R\colon T^t\text{ is {\rm not} strictly ergodic}\}
 \]
   is a countable set, cf. \cite[4.24.2]{G}.
\end{proof}

 We now specialize Theorem~\ref{mainsection2} to the commutative setting. If $(X,S)$ is a dynamical system, the topological entropy of $S$, which we will denote by $\mathrm{h}_{\mathrm{top}}(S)$, is either a nonnegative real number or infinite. Topological entropy is an invariant of conjugacy,  cf. \cite[Theorem~7.2]{W}.
 \begin{corollary}\label{topentropy}
 Let $(X,S)$ be a dynamical system and let $(Y,T)$ be the suspension flow associated to $(X,S)$. If $t\in\R$ then
 \[
     \topentr(T^t)=|t|\, \topentr ( S).
 \]
 \end{corollary}
 \begin{proof}
 Let $\alpha$ be the automorphism of $C(X)$ induced by $S$ and let $\phi$ be the action of $\R$ on $C(Y)$ induced by $T$. Then for each $t\in\R$, the Voiculescu-Brown entropy of $\phi_t$ and the topological entropy of $T^t$ coincide, cf.\ \cite[Proposition~4.8]{Vo}.
 By the same token, the Voiculescu-Brown entropy of $\alpha$ is equal to the topological entropy of $S$. Thus the corollary follows from Theorem~\ref{mainsection2}.
 \end{proof}
 
 We remark that Corollary~\ref{topentropy} can be obtained, without Theorem~\ref{mainsection2}, using tools from classical topological entropy, cf.\ \cite[Theorem~IV.5]{I-O_Thesis}. 
Using our results above we are able to give the following.
 
 \begin{example}\label{theexample}
     There are dynamical systems on (compact metric)  connected 1-dimensional spaces which are not flip conjugate yet their associated crossed products have isomorphic Elliott invariants.
 \end{example}
 
       Let $(X_1,S_1)$ be a strictly ergodic Cantor system and let $s\not =\topentr(S_1)$ be a nonnegative real number or infinity.
       By \cite[Theorem~6.1]{Su1} and \cite[Theorem~7.1]{Su2}, there is a  minimal Cantor system $(X_2,S_2)$ strongly orbit equivalent to $(X_1,S_2)$ and with entropy equal to $s$. Since orbit equivalence preserves the space of invariant Borel probability measures (cf.\ \cite[Theorem~2.2]{GPS}), we have that $(X_2,S_2)$ is also a strictly ergodic Cantor system.  Let $(Y_1,T_1)$ and $(Y_2,T_2)$ be the suspension flows of $(X_1,S_1)$ and $(X_2,S_2)$.  Then, by  Proposition~\ref{texists}, there is a $t\in\R$ such that $T_{1}^{t}$  and $T_{2}^{t}$ are both strictly ergodic.
  Hence $(Y_1,T_{1}^{t})$ and $(Y_2,T_{2}^{t})$ are dynamical systems on (compact metric) connected 1-dimensional spaces which are not flip conjugate (because they have different entropy, by Corollary~\ref{topentropy}) yet their associated crossed products have  isomorphic Elliott invariants (by Theorem~\ref{ElliottInv}).
  {\vspace{-.58cm}{\begin{flushright}$\Box$\end{flushright}}
 
\medskip
We conclude this paper with a brief discussion about the classifiability of the C*-algebras  whose Elliott invariants were calculated in Proposition~\ref{elliottinvariants}.\footnote{Note Added in Proof: Recently, H.\ Lin and N.\ C.\ Phillips have proved, in  \cite{LP}, a result  which shows that the C*-algebras of Theorem~\ref{ElliottInv} (and in particular, those of Example~\ref{theexample}) are in fact isomorphic.}
 We conjecture they are classifiable.  For the particular case of minimal rotations of suspensions arising from group rotations the result follows from a theorem of Gjerde and Johansen \cite{GJ}. Indeed, suppose that the dynamical systems $(X_1,S_1)$ and $(X_2,S_2)$ are minimal group rotations. Let $(Y_1,T_1)$ and $(Y_2,T_2)$ be the suspension flows associated to $(X_1,S_1)$ and $(X_2,S_2)$, respectively. Then it is easy to see that  $(Y_1,T_{1}^{t_1})$ and $(Y_2,T_{2}^{t_2})$ are also  group rotations for all $t_1,t_2\in\R$. Hence, when $t_1$ and $t_2$ are real numbers such that $T_{1}^{t_1}$ and $T_{2}^{t_2}$ are minimal (and so they are strictly ergodic, cf.\ \cite[Corollary~6.20]{W})  we may use  \cite[Theorem~6]{GJ} to obtain that   $C\sp\ast (\Z,Y_1,T_{1}^{t_1})$ and $C\sp\ast (\Z,Y_2,T_{2}^{t_2})$ are isomorphic if and only if $C\sp\ast (\Z,Y_1,T_{1}^{t_1})$ and $C\sp\ast (\Z,Y_2,T_{2}^{t_2})$ have isomorphic Elliott invariants. We have proved the following.

\begin{proposition}\label{classification}
Suppose that the dynamical systems $(X_1,S_1)$ and $(X_2,S_2)$ are minimal group rotations. Let $(Y_1,T_1)$ and $(Y_2,T_2)$ be the suspension flows associated to $(X_1,S_1)$ and $(X_2,S_2)$, respectively. If $t_1$ and $t_2$ are real numbers such that $T_{1}^{t_1}$ and $T_{2}^{t_2}$ are minimal then $C\sp\ast (\Z,Y_1,T_{1}^{t_1})$ and $C\sp\ast (\Z,Y_2,T_{2}^{t_2})$ are isomorphic if and only if $C\sp\ast (\Z,Y_1,T_{1}^{t_1})$ and $C\sp\ast (\Z,Y_2,T_{2}^{t_2})$ have isomorphic Elliott invariants.
 \end{proposition}
 
 Proposition~\ref{classification} does not cover the C*-algebras of Example~\ref{theexample}, unfortunately. However,  it does cover the important example of irrational rotation algebras:  An irrational rotation on $S^1$ can be thought as  a time-$t$ homeomorphism of the suspension built over the  one-point dynamical system. Since the dynamical system formed by the one-point space is a minimal  group rotation, Proposition~\ref{classification} applies.  Hence two irrational rotation algebras $A_{t_1}$ and $A_{t_2}$ are isomorphic if and only if they have isomorphic Elliott invariants. Since irrational rotations on the unit circle are strictly ergodic, we may use our Propostion~\ref{elliottinvariants} to obtain the Elliott invariants of $A_{t_1}$ and $A_{t_2}$. This invariants will be isomorphic if and only if $\Z + t_1\, \Z =\Z + t_2\, \Z$ if and only if $t_1$ have the same image as $\pm t_2$ in $\R/ \Z$.

\flushleft{\textit{\mbox{} \\
Department of Mathematics and Statistics,
University of Ottawa\\ 
585 King Edward Ave., Ottawa, Ontario, Canada K1N-6N5\\
E-mail: bitzaort@uottawa.ca\\
}}

\end{document}